\documentclass[graybox]{svmult}

\usepackage{tikz,cite,multirow,array,enumitem,diagbox,capt-of,hyperref,mathtools}
\usepackage{type1cm}        
\usepackage{makeidx}         
\usepackage{graphicx}        
\usepackage{multicol}        
\usepackage[bottom]{footmisc}

\usepackage{newtxtext}       %
\usepackage[varvw]{newtxmath}       

\newcommand\scalemath[2]{\scalebox{#1}{\mbox{\ensuremath{\displaystyle#2}}}}

\makeindex             

\begin{document}

\title*{Towards a Multigrid Preconditioner Interpretation of
  Hierarchical Poincaré-Steklov Solvers}
\titlerunning{A HPS Multigrid Preconditioner}
\author{Jos\'e Pablo Lucero Lorca \orcidID{0000-0002-9005-4146}}
\institute{\email{mail@pablo.world}}
%
%
\maketitle

\abstract{%
We revisit the Hierarchical Poincaré--Steklov (HPS) method in a
preconditioned iterative setting for variable-coefficient Helmholtz
problems with impedance boundary conditions. HPS is commonly presented
as a direct solver based on nested dissection and high-order
tensor-product discretizations; here we recast its hierarchical merge
tree as a multilevel preconditioner for the assembled skeleton
(trace) system. The main goal is to \emph{flexibilize} the final,
memory-intensive coarse stage of direct HPS by replacing the exact
coarse solve with a small, fixed amount of iterative work, thereby
exposing tunable trade-offs between memory footprint and time to
solution. Numerical experiments on a two-dimensional scattering
benchmark illustrate these trade-offs and compare against both
unpreconditioned GMRES and the classic direct HPS pipeline with an
exact coarse space.%
}

\section{Introduction}
The Hierarchical Poincaré--Steklov (HPS) method, introduced by
Martinsson~\cite{Martinsson2013,Martinsson2015}, is a direct solver
for elliptic boundary-value problems that combines nested dissection
with spectral element discretizations on tensor-product
grids. Subsequent work extended this framework to variable-coefficient
Helmholtz equations, showing high accuracy and efficiency at
scale~\cite{GillmanMartinsson2014,
  GillmanBarnettMartinsson2015,BabbGillmanHaoMartinsson2018}. In the
impedance-to-impedance (ItI) formulation---based on the discretization
of Despr\'es~\cite{Despres1991}---Dirichlet and Neumann traces are
replaced by local impedance maps, yielding a closed interface
representation well suited to high-frequency and heterogeneous media.

The present work places HPS in a preconditioned iterative setting,
where the hierarchical merging tree provides the multilevel
organization. Viewed equivalently as a nested-dissection solver for a
spectral element discretization, HPS naturally induces a multilevel
preconditioner, thereby unifying the direct and iterative
perspectives.

Related hierarchical and multilevel strategies for elliptic and
Helmholtz-type problems are numerous. On the multilevel side,
classical frameworks include methods such
as~\cite{CaiGoldsteinPasciak1993}. For hierarchical discretizations
and direct/hybridization-style solvers, composite spectral collocation
and related multidomain spectral approaches are surveyed in
\cite{GillmanBarnettMartinsson2015} and include, for example,
\cite{PfeifferKidderScheelTeukolsky2003}. For Helmholtz problems,
iterative approaches including shifted-Laplace preconditioning,
sweeping and domain decomposition methods, and multigrid variants are
reviewed in~\cite{GanderZhang2019}.

The present work targets preconditioning of the skeleton (trace)
system and leverages the modular HPS construction
of~\cite{OutrataLuceroLorca2025}. A key practical drawback of HPS
pipelines is their memory footprint and communication volume in
distributed-memory settings, due to the exchange of dense interface
operators (see~\cite{MeliaFortunatoGillmanONeil2025} and references
therein). Our aim is to \emph{add} flexibility by avoiding
construction of the most expensive coarse spaces: on the final level,
the exact solve is replaced by a small, fixed number of iterations of
an iterative method applied to the assembled last-level system. A
broader comparison with other iterative methods is deferred to future
work.

\section{Model problem}
We consider the variable-coefficient Helmholtz equation with impedance
boundary conditions
\begin{align}
  -\Delta u - \underbrace{\kappa^2 \left(1-b(\boldsymbol{x})\right)}_{:=c(\boldsymbol{x})}u = s(\boldsymbol{x}),\quad \boldsymbol{x}\in\Omega
  \quad \text{and} \quad
  \frac{\partial u}{\partial n} + i\eta u = t(\boldsymbol{x}),
  \quad \boldsymbol{x}\in\partial\Omega,
  \label{eqn:diffhelm}
\end{align}
where $\Omega=(0,1)^2\subset\mathbb{R}^2$ and $u:\Omega\to\mathbb{C}$
is the unknown field, $\eta\in\mathbb{R}$ chosen equal to
$\kappa\in\mathbb{R}$ the wavenumber, $b(\boldsymbol{x})$ a smooth
coefficient, and $s(\boldsymbol{x}), t(\boldsymbol{x})$ smooth source
and boundary data. Impedance boundary conditions of this form are
widely used in diffraction, acoustics, and electromagnetic scattering
\cite{KyurkchanSmirnova2016, KinslerFreyCoppensSanders2000,
  Lighthill1978, ChandlerWilde1997}; see also
\cite[\S1.1,\S1.2]{GrahamSauter2020} for an overview.

\section{Discretization}\label{sec:discretization}
Consider a structured spectral element mesh, \(\Omega=(0,1)^2\) is
divided into a square grid of square elements, each with a
tensor–product Gauss–Legendre–Lobatto (GLL) grid of order \(N\).  This
construction allows high-order local operators from the tensor product
of 1D differentiation and mass matrices while preserving continuity of
impedance data on shared edges (see \cite{Despres1991}). Local ItI
maps are assembled element-by-element and coupled through interface
conditions as described in the following sections (more detailed
expositions can be found in \cite{BeamsGillmanHewett2020,
  LuceroLorcaBeamsBeecroftGillman2024} and references therein).

\subsection{Local discretization}\label{sec:localdiscretization}

Each element problem is represented by
\begin{align}
  \tilde{L}
  =&
  K_x\otimes M_y + M_x\otimes K_y
   + \operatorname{diag}\!\bigl(c(x_i,y_j)\bigr)(M_x\otimes M_y),
\end{align}
where $w_i,w_j$ are GLL quadrature weights,
$M_x=\operatorname{diag}(w_i)$ and $M_y=\operatorname{diag}(w_j)$ are
1D GLL mass matrices, $K_x=D_x^{\top}M_xD_x$ and
\(K_y=D_y^{\top}M_yD_y\) are stiffness matrices, and \(D_x,D_y\) are
the 1D differentiation matrices.  The diagonal operator contains the
coefficient $c$ evaluated on the tensor grid
\(\{(x_i,y_j)\}_{i,j=1}^{N+1}\).

Following \cite{LuceroLorcaBeamsBeecroftGillman2024}, the corner nodes
are removed from the discretization, since they can be recalculated
later in post-processing --- this is a property of tensor-product
spectral methods. The boundary index sets are denoted
\(\iota_{l},\iota_{r},\iota_{b},\iota_{t}\) for the left, right,
bottom, and top edges, and their union is \(\iota_{\Gamma}\). The
inner index set, denoted \(\iota_{i}\), contains all remaining nodes
strictly inside the element. The outgoing and incoming impedance
operators are
$  \mathcal{I}_o=\left[\begin{smallmatrix}
    -D_x \otimes I\\
     D_x \otimes I\\
    -I \otimes D_y\\
     I \otimes D_y
   \end{smallmatrix}\right](\iota_{\Gamma},:)-\eta\,I(\iota_{\Gamma}~,~:~)$
 and
 $
    \mathcal{I}_i~=~\left[\begin{smallmatrix}
    -D_x \otimes I\\
     D_x \otimes I\\
    -I \otimes D_y\\
     I \otimes D_y
  \end{smallmatrix}\right]~(~\iota_{\Gamma}~,~:~)~+~\eta\,I(\iota_{\Gamma},:),$
where $I$ is the identity of appropriate size.

To apply incoming impedance conditions, the boundary rows of
$\tilde L$ are replaced by $\mathcal I_i$ to define $L$ as
$
  L(\iota_{\Gamma},:) := \mathcal I_i \text{ and } 
  L(\iota_i,:) := \tilde L(\iota_i,:) .
$
The local Impedance-to-Impedance operator and interior contribution are
\begin{align}
  T =& \mathcal I_o L^{-1} I(:,\iota_{\Gamma}) , &
  H =& \mathcal I_o L^{-1} I(:,\iota_i)\,\tilde b(\iota_i) ,
  \label{eqn:ItIOperator}
\end{align}
where $\tilde b$ contains the local right-hand-side values.

The operators $T$ and $H$ yield the closed impedance relation
$
\mathcal{I}_o~u(\iota_{\Gamma})~=~T~\mathcal{I}_i~u(\iota_{\Gamma})~+~H
$,
from which the full element solution follows by
\begin{align}
  L u =& b, 
\end{align}
where $b(\iota_{\Gamma}) = \mathcal I_i u(\iota_{\Gamma})$ and
$b(\iota_i) = \tilde b(\iota_i)$.

\subsection{Global discretization}\label{sec:globaldiscretization}
For each element $e$, let $\alpha,\beta\in\{l,r,b,t\}$ denote its sides.
The local relation between outgoing and incoming impedance data is
\begin{align}\label{eqn:localiti}
  \scalemath{0.8}{
  (\mathcal I_o u)^{(e)}_{\alpha}
  =
  \sum_{\beta\in\{l,r,b,t\}}
  T^{(e)}_{\alpha\beta}(\mathcal I_i u)^{(e)}_{\beta}
  + H^{(e)}_{\alpha}} ,
\end{align}
where $T^{(e)}_{\alpha\beta}\in\mathbb C^{(N-1)\times(N-1)}$ and
$H^{(e)}_{\alpha}\in\mathbb C^{(N-1)}$ are the operators defined in
eq. \eqref{eqn:ItIOperator}.

Let two elements $e_1$ and $e_2$ share an interior face. Transmission
conditions enforce continuity of impedance data across shared faces:
\begin{align}
  (\mathcal I_i u)^{(e_2)}_{\beta}
  =&
  (\mathcal I_o u)^{(e_1)}_{\alpha},
  &
  (\mathcal I_i u)^{(e_1)}_{\alpha}
  =&
  (\mathcal I_o u)^{(e_2)}_{\beta}.
\end{align}
Combining these with the local ItI maps gives the face system
\begin{align}
  \left[\begin{smallmatrix}
    I & -T^{(e_1)}_{\alpha\alpha}\\
    -T^{(e_2)}_{\beta\beta} & I
  \end{smallmatrix}\right]
  \left[\begin{smallmatrix}
    (\mathcal I_i u)^{(e_2)}_{\beta}\\
    (\mathcal I_i u)^{(e_1)}_{\alpha}
  \end{smallmatrix}\right]
  - 
  \left[\begin{smallmatrix}
    \sum_{\gamma\ne\alpha}
    T^{(e_1)}_{\alpha\gamma}(\mathcal I_i u)^{(e_1)}_{\gamma}\\
    \sum_{\gamma\ne\beta}
    T^{(e_2)}_{\beta\gamma}(\mathcal I_i u)^{(e_2)}_{\gamma}
  \end{smallmatrix}\right]
  =&
  \left[\begin{smallmatrix}
    H^{(e_1)}_{\alpha}\\
    H^{(e_2)}_{\beta}
  \end{smallmatrix}\right].
  \label{eqn:faceequations}
\end{align}

Assembling all face equations yields the sparse non-hermitian global
skeleton system \[M g = \mathrm{RHS},\] where $g$ collects all interior
incoming impedances and $\mathrm{RHS}$ stacks the local
$H^{(e)}_{\alpha}$ contributions. Physical boundary sides contribute
directly to the right-hand side.

\section{Solver}\label{sec:solver}

The HPS solver applies the nested-dissection procedure to the spectral
element system described above. This section details the face ordering
that enables its direct solution, later recasted as a relaxation
scheme.

\subsection{Nested dissection: Local scheme}

Let two elements $e_1$ and $e_2$ share an interior face $f$ through sides
$\alpha$ of $e_1$ and $\beta$ of $e_2$. Their face equations
(from~\eqref{eqn:faceequations}) are
\begin{align}
  \scalemath{0.8}{
  (\mathcal I_i u)^{(e_2)}_{\beta}
  - \sum_{\gamma\in\{l,r,b,t\}}
    T^{(e_1)}_{\alpha\gamma}(\mathcal I_i u)^{(e_1)}_{\gamma}
  =}&   \scalemath{0.8}{
      H^{(e_1)}_{\alpha}},\\
    \scalemath{0.8}{
  (\mathcal I_i u)^{(e_1)}_{\alpha}
  - \sum_{\gamma\in\{l,r,b,t\}}
    T^{(e_2)}_{\beta\gamma}(\mathcal I_i u)^{(e_2)}_{\gamma}
  =}&  \scalemath{0.8}{ H^{(e_2)}_{\beta}}.
\end{align}
We now group the internal and external incoming impedances as
$x = \left[\begin{smallmatrix}
    (\mathcal I_i u)^{(e_2)}_{\beta} &
    (\mathcal I_i u)^{(e_1)}_{\alpha}
  \end{smallmatrix}\right]^\top$
  and
  $y = \left[\begin{smallmatrix}
      (\mathcal I_i u)^{(e_1)}_{E_1} &
      (\mathcal I_i u)^{(e_2)}_{E_2}
    \end{smallmatrix}\right]^\top$
where $E_1=\{l,r,b,t\}\!\setminus\!\{\alpha\}$ and
$E_2=\{l,r,b,t\}\!\setminus\!\{\beta\}$.
With this notation the system becomes
\begin{align}
  \underbrace{
    \left[\begin{smallmatrix}
      I & -T^{(e_1)}_{\alpha\alpha}\\
      -T^{(e_2)}_{\beta\beta} & I
    \end{smallmatrix}\right]}_{\scriptstyle A}
  x
  =&
  \underbrace{
    \left[\begin{smallmatrix}
      T^{(e_1)}_{\alpha E_1} & 0\\
      0 & T^{(e_2)}_{\beta E_2}
    \end{smallmatrix}\right]}_{\scriptstyle B}
  y
  +
  \underbrace{
    \left[\begin{smallmatrix}
      H^{(e_1)}_{\alpha}\\
      H^{(e_2)}_{\beta}
    \end{smallmatrix}\right]}_{\scriptstyle h} .
  \label{eqn:localfacesystem}
\end{align}
Eliminating $x$ gives $x = A^{-1}B y + A^{-1}h$, substituting into
the outgoing relations \eqref{eqn:localiti} produces the fused pair
operator
\begin{align}
  \scalemath{0.8}{
T_{\mathrm{pair}}
  =}&\scalemath{0.8}{
  \underbrace{D}_{\scriptstyle
    \left[\begin{smallmatrix}
      T^{(e_1)}_{E_1E_1} & 0\\
      0 & T^{(e_2)}_{E_2E_2}
    \end{smallmatrix}\right]}
  -
  \underbrace{C}_{\scriptstyle
    \left[\begin{smallmatrix}
      T^{(e_1)}_{E_1\alpha} & 0\\
      0 & T^{(e_2)}_{E_2\beta}
    \end{smallmatrix}\right]}
  A^{-1} B},
  &\qquad
    \scalemath{0.8}{
  H_{\mathrm{pair}}
    =} &
         \scalemath{0.8}{
  \left[\begin{smallmatrix}
    H^{(e_1)}_{E_1}\\
    H^{(e_2)}_{E_2}
  \end{smallmatrix}\right]
  +
  C A^{-1}
  \left[\begin{smallmatrix}
    H^{(e_1)}_{\alpha}\\
    H^{(e_2)}_{\beta}
        \end{smallmatrix}\right]}.
     \label{eqn:globalfacesystem}
\end{align}
where $T_{\mathrm{pair}}$ is clearly a Schur complement.

\subsection{Nested dissection: Global scheme}\label{sec:NDGlobal}

\begin{figure}
  \[
    \begin{array}{|c|}
        \hline
      \begin{array}{cc}
      \begin{array}{c}
        \text{\bf Grid 1} \\
    \resizebox{0.28\textwidth}{!}{%
    \begin{tikzpicture}
      \draw[step=1cm,thick] (0,0) grid (4,4);
      \draw[thick,gray!50] (0,0) rectangle (4,4);
      \tikzset{
        face/.style={circle,draw,fill=white,minimum size=9pt,inner
          sep=0pt,font=\tiny}
      }
      \newcommand{\face}[3]{%
        \pgfmathtruncatemacro{\ia}{mod(#2-1,4)}%
        \pgfmathtruncatemacro{\ja}{(#2-1)/4}%
        \pgfmathtruncatemacro{\ib}{mod(#3-1,4)}%
        \pgfmathtruncatemacro{\jb}{(#3-1)/4}%
        \pgfmathsetmacro{\x}{(\ia+\ib)/2+0.5}%
        \pgfmathsetmacro{\y}{(\ja+\jb)/2+0.5}%
        \node[face] at (\x,\y){#1};}
      \face{1}{1}{2}\face{2}{3}{4}\face{3}{5}{6}\face{4}{7}{8}
      \face{5}{9}{10}\face{6}{11}{12}\face{7}{13}{14}\face{8}{15}{16}
      \face{9}{1}{5}\face{10}{2}{6}\face{11}{9}{13}\face{12}{10}{14}
      \face{13}{3}{7}\face{14}{4}{8}\face{15}{11}{15}\face{16}{12}{16}
      \face{17}{2}{3}\face{18}{6}{7}\face{19}{10}{11}\face{20}{14}{15}
      \face{21}{5}{9}\face{22}{6}{10}\face{23}{7}{11}\face{24}{8}{12}
    \end{tikzpicture}
       }
   \end{array}
  &
    \begin{array}{cc}
    \begin{array}{c}
      M_1 = \left(\begin{matrix}A_1&B_1\\C_1&D_1\end{matrix}\right) =
    \end{array}
    \end{array}
    \begin{array}{c}
      \left(
      \begin{array}{cc}
        \begin{tabular}{c}
          \resizebox{0.07\textwidth}{!}{%
          \begin{tikzpicture}[x=0.4cm,y=0.4cm]
            \foreach \i/\j in {1/1,2/2,3/3,4/4,5/5,6/6,7/7,8/8}{
              \fill (\j-1,{8-\i}) rectangle (\j,{9-\i});
            }
            \draw[step=1,gray!30,thick] (0,0) grid (8,8);
          \end{tikzpicture}}
        \end{tabular}
      &
        \begin{tabular}{c}
          \resizebox{0.14\textwidth}{!}{%
          \begin{tikzpicture}[x=0.4cm,y=0.4cm]
          \foreach \i/\j in {
            1/1,1/2,1/9,
            2/5,2/6,2/9,
            3/1,3/2,3/10,3/13,3/14,
            4/5,4/6,4/10,4/15,4/16,
            5/3,5/4,5/11,5/13,5/14,
            6/7,6/8,6/11,6/15,6/16,
            7/3,7/4,7/12,
            8/7,8/8,8/12
          }{
            \fill (\j-1,{8-\i}) rectangle (\j,{9-\i});
          }
            \draw[step=1,gray!30,thick] (0,0) grid (16,8);
        \end{tikzpicture}}
      \end{tabular}
      \\
      \begin{tabular}{c}
        \resizebox{0.07\textwidth}{!}{%
        \begin{tikzpicture}[x=0.4cm,y=0.4cm]
          \foreach \i/\j in {
            1/1,1/3,
            2/1,2/3,
            3/5,3/7,
            4/5,4/7,
            5/2,5/4,
            6/2,6/4,
            7/6,7/8,
            8/6,8/8,
            9/1,9/2,
            10/3,10/4,
            11/5,11/6,
            12/7,12/8,
            13/3,13/5,
            14/3,14/5,
            15/4,15/6,
            16/4,16/6
          }{
            \fill (\j-1,{16-\i}) rectangle (\j,{17-\i});
          }
          \draw[step=1,gray!30,thick] (0,0) grid (8,16);
        \end{tikzpicture}}
      \end{tabular}
      &
      \begin{tabular}{c}
        \resizebox{0.14\textwidth}{!}{%
        \begin{tikzpicture}[x=0.4cm,y=0.4cm]
          \foreach \i/\j in {
            1/1,1/13,
            2/2,2/9,2/10,2/14,
            3/3,3/13,
            4/4,4/11,4/12,4/14,
            5/5,5/9,5/10,5/15,
            6/6,6/16,
            7/7,7/11,7/12,7/15,
            8/8,8/16,
            9/2,9/5,9/9,
            10/2,10/5,10/10,10/14,10/15,
            11/4,11/7,11/11,11/14,11/15,
            12/4,12/7,12/12,
            13/1,13/3,13/13,
            14/2,14/4,14/10,14/11,14/14,
            15/5,15/7,15/10,15/11,15/15,
            16/6,16/8,16/16
          }{
            \fill (\j-1,{16-\i}) rectangle (\j,{17-\i});
          }
          \draw[step=1,gray!30,thick] (0,0) grid (16,16);
        \end{tikzpicture}}
      \end{tabular}
    \end{array}
    \right)
   \end{array}
  \end{array}
\\\hline
\begin{array}{c|c|c}
  \begin{array}{c}
    \begin{array}{c}
      \text{\bf Grid 2} \\
      \resizebox{0.28\textwidth}{!}{%
      \begin{tikzpicture}
        \draw[thick,gray!50] (0,0) rectangle (4,4);
        \foreach \y in {1,2,3} \draw[thick] (0,\y) -- (4,\y);
        \draw[thick] (2,0) -- (2,4);
        \tikzset{
          face/.style={circle,draw,fill=white,minimum size=9pt,inner
            sep=0pt,font=\tiny}
        }
        \newcommand{\face}[3]{%
          \pgfmathtruncatemacro{\ia}{mod(#2-1,4)}%
          \pgfmathtruncatemacro{\ja}{(#2-1)/4}%
          \pgfmathtruncatemacro{\ib}{mod(#3-1,4)}%
          \pgfmathtruncatemacro{\jb}{(#3-1)/4}%
          \pgfmathsetmacro{\x}{(\ia+\ib)/2+0.5}%
          \pgfmathsetmacro{\y}{(\ja+\jb)/2+0.5}%
          \node[face] at (\x,\y){#1};}
        \face{1}{1}{5}\face{2}{2}{6}\face{3}{9}{13}\face{4}{10}{14}
        \face{5}{3}{7}\face{6}{4}{8}\face{7}{11}{15}\face{8}{12}{16}
        \face{9}{2}{3}\face{10}{6}{7}\face{11}{10}{11}\face{12}{14}{15}
        \face{13}{5}{9}\face{14}{6}{10}\face{15}{7}{11}\face{16}{8}{12}
      \end{tikzpicture}
      }
    \end{array}
    \\
    \begin{array}{c}
      M_2 = D_1 - C_1 A_1^{-1} B_1
    \end{array}
    \\
    = \left(
    \begin{array}{c}
      \begin{array}{cc}
        \begin{tabular}{c}
          \resizebox{0.07\textwidth}{!}{%
          \begin{tikzpicture}[x=0.4cm,y=0.4cm]
            \foreach \i/\j in {
              1/1,1/2,2/1,2/2,
              3/3,3/4,4/3,4/4,
              5/5,5/6,6/5,6/6,
              7/7,7/8,8/7,8/8
            }{
              \fill (\j-1,{8-\i}) rectangle (\j,{9-\i});
            }
            \draw[step=1,gray!30,thick] (0,0) grid (8,8);
          \end{tikzpicture}}
        \end{tabular}
        &
          \begin{tabular}{c}
            \resizebox{0.07\textwidth}{!}{%
            \begin{tikzpicture}[x=0.4cm,y=0.4cm]
              \foreach \i/\j in {
                1/1,1/2,2/1,2/2,
                1/5,1/6,2/5,2/6,
                3/3,3/4,4/3,4/4,
                3/5,3/6,4/5,4/6,
                5/1,5/2,6/1,6/2,
                5/7,5/8,6/7,6/8,
                7/3,7/4,8/3,8/4,
                7/7,7/8,8/7,8/8
              }{
                \fill (\j-1,{8-\i}) rectangle (\j,{9-\i});
              }
              \draw[step=1,gray!30,thick] (0,0) grid (8,8);
            \end{tikzpicture}}
          \end{tabular}
        \\
        \begin{tabular}{c}
          \resizebox{0.07\textwidth}{!}{%
          \begin{tikzpicture}[x=0.4cm,y=0.4cm]
            \foreach \i/\j in {
              1/1,1/2,2/1,2/2,
              1/5,1/6,2/5,2/6,
              3/3,3/4,4/3,4/4,
              3/7,3/8,4/7,4/8,
              5/1,5/2,6/1,6/2,
              5/3,5/4,6/3,6/4,
              7/5,7/6,8/5,8/6,
              7/7,7/8,8/7,8/8
            }{
              \fill (\j-1,{8-\i}) rectangle (\j,{9-\i});
            }
            \draw[step=1,gray!30,thick] (0,0) grid (8,8);
          \end{tikzpicture}}
        \end{tabular}
        &
          \begin{tabular}{c}
            \resizebox{0.07\textwidth}{!}{%
            \begin{tikzpicture}[x=0.4cm,y=0.4cm]
              \foreach \i/\j in {
                1/1,2/2,3/3,4/4,
                2/5,2/6,2/7,2/8,
                3/5,3/6,3/7,3/8,
                5/2,5/3,
                6/2,6/3,
                7/2,7/3,
                8/2,8/3,
                5/5,5/6,
                6/5,6/6,
                7/7,7/8,
                8/7,8/8
              }{
                \fill (\j-1,{8-\i}) rectangle (\j,{9-\i});
              }
              \draw[step=1,gray!30,thick] (0,0) grid (8,8);
            \end{tikzpicture}}
          \end{tabular}
      \end{array}
    \end{array}
    \right)
  \end{array}
        &
          \begin{array}{c}
            \text{\bf Grid 3} \\
    \begin{tabular}{c}
      \resizebox{0.28\textwidth}{!}{%
      \begin{tikzpicture}
        \draw[thick,gray!50] (0,0) rectangle (4,4);
        \draw[thick] (0,2) -- (4,2);
        \draw[thick] (2,0) -- (2,4);
        \tikzset{
          face/.style={circle,draw,fill=white,minimum size=9pt,inner
            sep=0pt,font=\tiny}
        }
        \newcommand{\face}[3]{%
          \pgfmathtruncatemacro{\ia}{mod(#2-1,4)}%
          \pgfmathtruncatemacro{\ja}{(#2-1)/4}%
          \pgfmathtruncatemacro{\ib}{mod(#3-1,4)}%
          \pgfmathtruncatemacro{\jb}{(#3-1)/4}%
          \pgfmathsetmacro{\x}{(\ia+\ib)/2+0.5}%
          \pgfmathsetmacro{\y}{(\ja+\jb)/2+0.5}%
          \node[face] at (\x,\y){#1};}
        \face{1}{2}{3}
        \face{2}{6}{7}
        \face{3}{10}{11}
        \face{4}{14}{15}
        \face{5}{5}{9}
        \face{6}{6}{10}
        \face{7}{7}{11}
        \face{8}{8}{12}
      \end{tikzpicture}
      }
    \end{tabular}
    \\
    \begin{array}{c}
      M_3 = D_2 - C_2 A_2^{-1} B_2
    \end{array}
    \\
    =\left(
    \begin{array}{c}
      \begin{tabular}{cc}
        \begin{tabular}{c}
          \resizebox{0.07\textwidth}{!}{%
          \begin{tikzpicture}[x=0.4cm,y=0.4cm]
            \foreach \i/\j in {
              1/1,2/2
            }{
              \fill (\j-1,{2-\i}) rectangle (\j,{3-\i});
            }
            \draw[step=1,gray!30,thick] (0,0) grid (2,2);
          \end{tikzpicture}}
        \end{tabular}
        &
          \begin{tabular}{c}
            \resizebox{0.07\textwidth}{!}{%
            \begin{tikzpicture}[x=0.4cm,y=0.4cm]
              \foreach \i/\j in {
                1/1,1/2,2/1,2/2
              }{
                \fill (\j-1,{2-\i}) rectangle (\j,{3-\i});
              }
              \draw[step=1,gray!30,thick] (0,0) grid (2,2);
            \end{tikzpicture}}
          \end{tabular}
        \\
        \begin{tabular}{c}
          \resizebox{0.07\textwidth}{!}{%
          \begin{tikzpicture}[x=0.4cm,y=0.4cm]
              \foreach \i/\j in {
                1/1,1/2,2/1,2/2
              }{
                \fill (\j-1,{2-\i}) rectangle (\j,{3-\i});
              }
              \draw[step=1,gray!30,thick] (0,0) grid (2,2);
          \end{tikzpicture}}
        \end{tabular}
        &
        \begin{tabular}{c}
          \resizebox{0.07\textwidth}{!}{%
          \begin{tikzpicture}[x=0.4cm,y=0.4cm]
            \foreach \i/\j in {
              1/1,2/2
            }{
              \fill (\j-1,{2-\i}) rectangle (\j,{3-\i});
            }
            \draw[step=1,gray!30,thick] (0,0) grid (2,2);
          \end{tikzpicture}}
        \end{tabular}
      \end{tabular}
    \end{array}
    \right)
  \end{array}
  &
  \begin{array}{c}
    \text{\bf Grid 4} \\
    \begin{tabular}{c}
      \resizebox{0.28\textwidth}{!}{%
      \begin{tikzpicture}
        \draw[thick,gray!50] (0,0) rectangle (4,4);
        \draw[thick] (0,2) -- (4,2);
        \tikzset{
          face/.style={circle,draw,fill=white,minimum size=9pt,inner
            sep=0pt,font=\tiny}
        }
        \newcommand{\face}[3]{%
          \pgfmathtruncatemacro{\ia}{mod(#2-1,4)}%
          \pgfmathtruncatemacro{\ja}{(#2-1)/4}%
          \pgfmathtruncatemacro{\ib}{mod(#3-1,4)}%
          \pgfmathtruncatemacro{\jb}{(#3-1)/4}%
          \pgfmathsetmacro{\x}{(\ia+\ib)/2+0.5}%
          \pgfmathsetmacro{\y}{(\ja+\jb)/2+0.5}%
          \node[face] at (\x,\y){#1};}
        \face{1}{5}{9}
        \face{2}{6}{10}
        \face{3}{7}{11}
        \face{4}{8}{12}
      \end{tikzpicture}
      }
    \end{tabular}
    \\
    \begin{array}{c}
      M_4 = D_3 - C_3 A_3^{-1} B_3
    \end{array}
    \\
    =\left(
    \begin{array}{c}
      \begin{tabular}{cc}
          \begin{tabular}{c}
            \resizebox{0.07\textwidth}{!}{%
            \begin{tikzpicture}[x=0.4cm,y=0.4cm]
              \foreach \i/\j in {
                1/1,1/2,2/1,2/2
              }{
                \fill (\j-1,{2-\i}) rectangle (\j,{3-\i});
              }
              \draw[step=1,gray!30,thick] (0,0) grid (2,2);
            \end{tikzpicture}}
          \end{tabular}
        &
          \begin{tabular}{c}
            \resizebox{0.07\textwidth}{!}{%
            \begin{tikzpicture}[x=0.4cm,y=0.4cm]
              \foreach \i/\j in {
                1/1,1/2,2/1,2/2
              }{
                \fill (\j-1,{2-\i}) rectangle (\j,{3-\i});
              }
              \draw[step=1,gray!30,thick] (0,0) grid (2,2);
            \end{tikzpicture}}
          \end{tabular}
        \\
        \begin{tabular}{c}
          \resizebox{0.07\textwidth}{!}{%
          \begin{tikzpicture}[x=0.4cm,y=0.4cm]
              \foreach \i/\j in {
                1/1,1/2,2/1,2/2
              }{
                \fill (\j-1,{2-\i}) rectangle (\j,{3-\i});
              }
              \draw[step=1,gray!30,thick] (0,0) grid (2,2);
          \end{tikzpicture}}
        \end{tabular}
        &
          \begin{tabular}{c}
            \resizebox{0.07\textwidth}{!}{%
            \begin{tikzpicture}[x=0.4cm,y=0.4cm]
              \foreach \i/\j in {
                1/1,1/2,2/1,2/2
              }{
                \fill (\j-1,{2-\i}) rectangle (\j,{3-\i});
              }
              \draw[step=1,gray!30,thick] (0,0) grid (2,2);
            \end{tikzpicture}}
          \end{tabular}
      \end{tabular}
    \end{array}
    \right)
  \end{array}
\end{array}
\\\hline
\end{array}
\]
\caption{Face merging and sparsity patterns for a $4\times 4$ element
  mesh
  {\bf Grid 1:} Faces 1 to 8 are eliminated, merging pairs of elements. These
    faces' dofs form the top left $1\times 1$-face-block diagonal part
    of $M_1$ since they are not linked \emph{directly} between each
    other, but through another face, e.g. face 1 is related to face 2
    through face 17.
  {\bf Grid 2:}Faces 1 to 8 are eliminated by pairs, merging $1\times 2$
    subdomains by one of their largest sides. These faces' dofs form the
    top left $2\times 2$-face-block diagonal matrix.
  {\bf Grid 3:}Faces 1 to 4 are eliminated by pairs, merging $2\times 2$
    subdomains. Thes faces' dofs form the top left
    $2\times 2$-face-block diagonal matrix.
  {\bf Grid 4:}Faces 1 to 4 are now fully coupled, $M_4$ is dense.
  }
\label{fig:sparsitypatterns}
\end{figure}
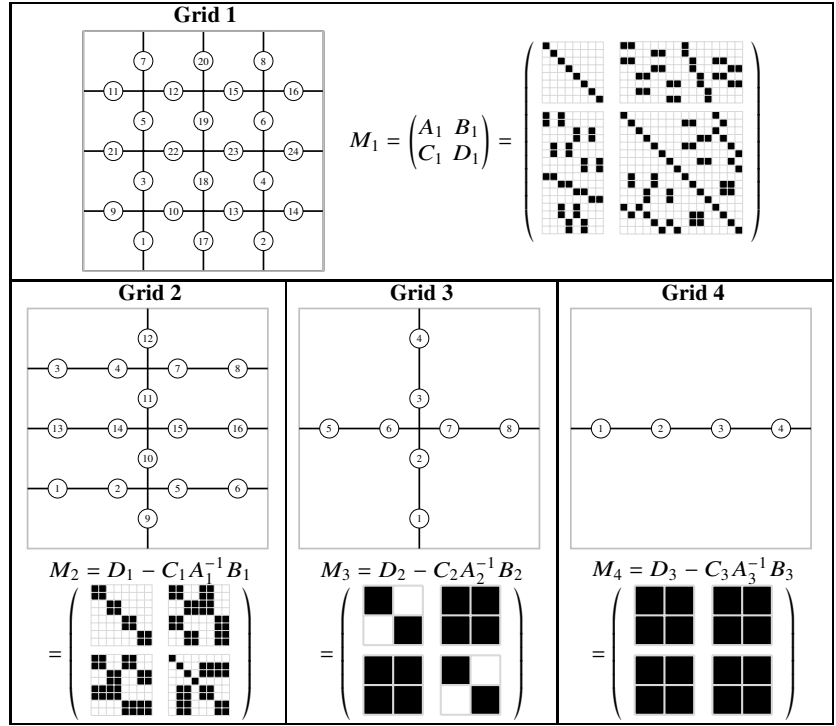

Figure \ref{fig:sparsitypatterns}
illustrates the face-merging procedure for the skeleton system $M_1$
on a $4\times 4$ element mesh.  The sequence Grid~1--Grid~4 is nested
dissection in reverse.  It lists, from fine to coarse, the elimination
sets used by the solver.  The key consequence of using this hierarchy
is that at every level $\ell$ the eliminated-face block $A_\ell$ is
block diagonal, with blocks growing as subdomains are merged.
Applying $A_\ell^{-1}$ to form the next Schur complement has
controlled cost, while the fill is pushed to coarser levels.

Let $f_i$ denote the face labeled $i$ in Grid~1.  The first
elimination set is $\mathcal F_1=\{f_1,\dots,f_8\}$.  It is maximal
among the nested dissection ordered subsets of Grid~1 faces that are
pairwise element-disjoint.  Thus the faces in $\mathcal F_1$ do not
couple directly and can be eliminated independently.  Eliminating
$\mathcal F_1$ merges the element pairs adjacent to these faces,
producing the $(1\times 2)$-element subdomains shown in Grid~2. Given
the nested dissection face ordering, $M_1$ is partitioned as
$M_1=\begin{psmallmatrix}A_1&B_1\\ C_1&D_1\end{psmallmatrix}$ and
eliminating $\mathcal F_1$ yields $M_2=D_1-C_1A_1^{-1}B_1$, the
reduced face system shown in Grid~2.  In this example the diagonal
blocks of $M_2$ are the fused ItI operators for the merged element
pairs.  The remaining nonzeros in $M_2$ encode couplings between these
merged pairs across the faces of Grid~2.  The coarse grids in
Grid~2--Grid~4 are tensor-product in the geometric sense that the
merges produce axis-aligned rectangular subdomains arranged in a
Cartesian tiling.

At level $\ell\ge 1$ we analogously eliminate a maximal
element-disjoint set $\mathcal F_\ell$ in the level-$\ell$ face graph
and obtain $M_{\ell+1}=D_\ell-C_\ell A_\ell^{-1}B_\ell$.  The
eliminated unknowns at level $\ell$ decompose into disjoint
merged-subdomain groups.  Accordingly $A_\ell$ is block diagonal, with
one block per merged subdomain at that level.  Grid~3 and Grid~4
visualize the next two levels.  The blocks in $A_\ell$ grow, and the
Schur complements become progressively denser as the hierarchy
coarsens.  Repeating this construction and using the associativity of
Schur complements \cite{CrabtreeHaynsworth1969} yields the
nested-dissection factorization of the spectral-element face system.

\subsection{Solver recast as a multigrid relaxation scheme}\label{sec:relaxationscheme}

The block-inverse relation introduced
in~\cite{LuceroLorcaRosenbergJankovMcCoidGander2025} takes the form
\begin{align}
  \scalemath{0.9}{
  M^{-1} = 
  \begin{bmatrix}
    A & B\\
    C & D
  \end{bmatrix}^{-1}
  =
  \underbrace{\begin{bmatrix}
    A^{-1} & 0\\
    0 & 0
  \end{bmatrix}}_{F}
  +
  \underbrace{\begin{bmatrix}
    -A^{-1}B\\
    I
  \end{bmatrix}}_{P}
  \underbrace{\left(D - C A^{-1} B\right)^{-1}}_{S^{-1}}
  \underbrace{\begin{bmatrix}
    0 & I
  \end{bmatrix}}_{R}
  \bigg(
    I - M
    \underbrace{\begin{bmatrix}
      A^{-1} & 0\\
      0 & 0
    \end{bmatrix}}_{F}
          \bigg).
          }
\end{align}
This identity motivates the definition of a
recursive multigrid algorithm without post-smoothing rather than a
single relaxation step: the local inversion $A^{-1}$ acts as a
smoother, and the reduced system $S$ defines the next level.  The
recursive iteration $\mathrm{MG}(M)= F + P S^{-1} R (I - M F)$,
where $S^{-1}$ is obtained by applying the same procedure to $S$. A
single coarse call yields a V-cycle; multiple ones define a
$\gamma$-cycle—both fully consistent with the hierarchical merging in
the HPS method. We employ MG as a preconditioner for flexible GMRES,
with the coarse solve performed by a fixed number of unpreconditioned
GMRES iterations.

\section{Numerical experiments}

We consider one of the problems
from~\cite{LuceroLorcaBeamsBeecroftGillman2024}, with
$b(\mathbf{x}) = 1.5e^{-160[(x-0.5)^2+(y-0.5)^2]}$ and
$s(\mathbf{x}) = -\kappa^2 b(\mathbf{x}) e^{i\kappa x}$, representing
scattering by a Gaussian bump.  We use polynomial degree~16, a
residual tolerance of~$10^{-8}$, and a frequency giving~9.6~points per
wavelength, yielding about~$10^{-7}$ accuracy (verifying the estimate
in \cite{LuceroLorcaBeamsBeecroftGillman2024}) for roughly one million
degrees of freedom before skeletonization.

Figure~\ref{fig:solution} shows the solution, and
Table~\ref{tab:results} reports results obtained in \textsc{MATLAB},
varying the number of levels. The table lists memory footprint, build
time, total iterations, and solve time for different fixed coarse
iteration counts and $\gamma$ values. The face sets used to build the
multilevel grids are those described in
Section~\ref{sec:NDGlobal}. The problem was run on a laptop with~32 GB
RAM and a hybrid processor (6 hyper-threading cores @ 4.7 GHz and 8
cores @ 3.5 GHz). Although cache effects favor certain configurations,
an overall timing trend can be observed. The method demonstrates that
performance can be tuned to available memory and the number of solves
required, while being faster than the unpreconditioned case in many
configurations. For context, we compare against the classic direct HPS
method, which uses an exact coarse space. It can be observed that
accepting a few iterations can save a significant amount of memory
footprint.

\section{Conclusion}
We provide a flexible iterative variant of the otherwise direct HPS
method for variable-coefficient Helmholtz problems arising, for
instance, in wave propagation and geological prospection: on different
coarse levels we replace the exact solve by a small, fixed
number of FGMRES iterations, explicitly trading a few Krylov steps for
a reduced coarse-space memory footprint.  The large-scale (including
3D) assessments in~\cite{LuceroLorcaBeamsBeecroftGillman2024,
  MeliaFortunatoGillmanONeil2025} indicate that the coarse-level
operators can dominate memory, even when factorized/compressed via
SVD-type techniques; our goal is to mitigate this bottleneck.  Since
current architectures tend to be more compute-rich than memory-rich,
we expect this tradeoff to enable larger problem sizes at fixed
accuracy by paying only a few additional FGMRES iterations.

\begin{table}[h!]
\centering
\begin{minipage}[t]{0.53\textwidth}
  \vspace{-3cm}
  \caption{{\bf PMem:} Preconditioner Memory Footprint [MB], {\bf It:}
    Flexible GMRES iterations with restart at 60, {\bf Bt:} Build time [s],
    {\bf St:} Solve time [s], {\bf c.i.:} coarse GMRES
    iterations. Results for $10^6$ dofs at $9.6$ points per
    wavelength.}
\label{tab:results}
\vspace{0cm}
\begin{tabular}{|l|c|c|c|c|}
\hline
Case & PMem & It & Bt & St \\ \hline
Unpreconditioned   & 0    & 669 & 0  & 85 \\ \hline
Exact coarse space & 3108 & 1   & 75 & 4  \\ \hline
\end{tabular}
\par\vspace{1mm}
\begin{tabular}{|r|r|r|rr|rr|rr|rr|rr|rr|}
\cline{4-15}
\multicolumn{3}{c}{} & \multicolumn{6}{|c|}{$\gamma = 1$} & \multicolumn{6}{c|}{$\gamma = 2$} \\ \hline
\multirow{2}{*}{\#levels} & \multirow{2}{*}{PMem} & \multirow{2}{*}{Bt} &
\multicolumn{2}{c|}{4 c.i.} & \multicolumn{2}{c|}{5 c.i.} &
\multicolumn{2}{c|}{6 c.i.} & \multicolumn{2}{c|}{2 c.i.} &
\multicolumn{2}{c|}{3 c.i.} & \multicolumn{2}{c|}{4 c.i.} \\ \cline{4-15}
 & & & It & St & It & St & It & St & It & St & It & St &
 It & St \\ \hline\hline
2  & 46   &  6 & 37 & 53 & 22 & 44 & 16 & 45 & 83 & 71 & 32 & 44 & 18 & 39 \\ \hline
3  & 460  & 15 & 23 & 42 & 15 & 40 & 11 & 40 & 24 & 55 & 11 & 42 &  7 & 42 \\ \hline
4  & 805  & 20 & 18 & 30 & 12 & 27 &  9 & 27 & 11 & 48 &  6 & 41 &  4 & 43 \\ \hline
5  & 1202 & 27 & 13 & 36 &  9 & 34 &  7 & 36 &  5 & 55 &  \textbf{3} & \textbf{71} &  1 & 52 \\ \hline
6  & 1527 & 31 & 11 & 22 &  7 & 18 &  5 & 16 &  2 & 47 &  1 & 46 & \multicolumn{2}{c}{} \\ \cline{1-13}
7  & 1897 & 38 &  9 & 28 &  6 & 26 &  4 & 23 &  \textbf{1} & \textbf{90} & \multicolumn{4}{c}{} \\ \cline{1-11}
8  & 2185 & 43 &  8 & 19 &  5 & 15 &  4 & 15 & \multicolumn{6}{c}{} \\ \cline{1-9}
9  & 2502 & 45 &  8 & 28 &  5 & 23 &  4 & 24 & \multicolumn{6}{c}{} \\ \cline{1-9}
10 & 2724 & 52 &  7 & 19 &  5 & 17 &  4 & 15 & \multicolumn{6}{c}{} \\ \cline{1-9}
11 & 2946 & 63 &  \textbf{7} & \textbf{26} &  \textbf{5} & \textbf{24} &  \textbf{4} & \textbf{24} & \multicolumn{6}{c}{} \\ \cline{1-9}
12 & 3051 & 67 &  3 & 11 &  2 &  8 &  1 &  6 & \multicolumn{6}{c}{} \\ \cline{1-9}
\end{tabular}
\vspace{-5mm}
\end{minipage}%
\hfill
\begin{minipage}[t]{0.35\textwidth}
  \centering
  \includegraphics[width=0.8\textwidth]{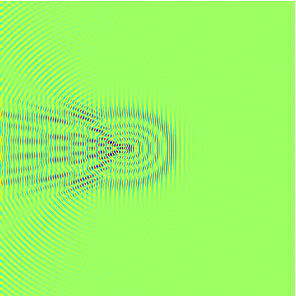}
  \captionof{figure}{Solution of the variable-coefficient Helmholtz problem.}
  \label{fig:solution}

  \vspace{2mm} {\small
    \noindent\textbf{Notes.}\par
    \raggedright $\eta$ is chosen equal to $\kappa$ to obtain
    $10^6$ dofs while yielding $\approx 10^{-7}$ accuracy; the
    estimate in~\cite{LuceroLorcaBeamsBeecroftGillman2024} was
    confirmed in these runs.  FGMRES is initialized with the zero
    vector. PMem includes the full memory usage of the
    program. Restart~60 was selected to avoid exhausting laptop
    memory. Bold entries satisfy $B_t+S_t>85$. \par }
\end{minipage}
\end{table}

\bibliographystyle{spmpsci} 
\bibliography{refs}
\end{document}